\documentclass[12pt,a4paper]{article}
\usepackage{amsmath,amssymb}
\usepackage{bm}

\makeatletter
\newcommand{\mysloppy}{\tolerance 9999 \hfuzz .5\p@ \vfuzz .5\p@}
\makeatother

\newtheorem{thm}{Theorem}[section]
\newtheorem{fact}[thm]{Fact}

\newtheorem{example}[thm]{Example}

\def\define{\mathrel{:=}}

\newcounter{myequation}[section]

\newcommand{\gor}{Gorenstein}

\newcommand{\covered}{\mathrel{<\!\!\!\cdot}}
\newcommand{\meet}{\wedge}
\newcommand{\join}{\vee}

\bmdefine{\NNN}{N}
\begin{document}\mysloppy
%\vspace*{2cm}
\begin{center}\LARGE
On the generating poset of Schubert cycles and
the characterization of
\gor\ property
\end{center}
\begin{center}\large
Mitsuhiro 
{\sc
Miyazaki}%
%\\
%\normalsize
\footnote{Dept.\ Math.\ Kyoto University of Education,
\
1 Fukakusa-Fujinomori-cho, Fushimi-ku, Kyoto,
612-8522 Japan,
\
E-mail:\tt
g53448@kyokyo-u.ac.jp}
\end{center}

\begin{trivlist}\item[]\small
{\bf Abstract:}\quad
The homogeneous coordinate ring of a Schubert variety
(a Schubert cycle for short) is an 
%ASL 
algebra with straightening law
generated by a
distributive lattice.
This paper gives a simple method to study the set of all the 
join-irreducible elements of this distributive lattice,
and gives a simple proof of the 
criterion 
%characterization
of the \gor{}
property of Schubert cycles.
\end{trivlist}
\begin{trivlist}\item[]\small
\gor{} property
\end{trivlist}

\section{Introduction}
In this note, all rings and algebras are assumed to be
commutative with identity element.
Let $m$ and $n$ be integers with $1\leq m\leq n$,
$B$ a ring and $S$ a $B$-algebra.
For an $m\times n$ matrix $M$ with entries in $S$,
we denote by $G(M)$ the $B$-subalgebra of $S$
generated by all the maximal minors of $M$.

From now on, we assume that $B$ is a field.
Let $X$ be the $m\times n$ matrix of indeterminates.
Then it is well known that $G(X)$ is the homogeneous
coordinate ring of the Grassmann variety $G_{m,n}$
of $m$-planes in $n$-space (\cite{hp}, \cite{dep2}, \cite{bv}).
And it is also well-known that $G(X)$ is a graded 
algebra with straightening law (graded ASL for short) over $B$
generated by the poset $\Gamma(X)$, where
\[
\Gamma(X)\define\{[a_1,\ldots,a_m]\mid
1\leq a_1<\cdots< a_m\leq n\}
\]
(\cite{dep2}, \cite{bv}).
The order of $\Gamma(X)$ is defined by
\[
[a_1,\ldots,a_m]\leq[b_1,\ldots,b_m]
\stackrel{\rm def}{\Longleftrightarrow}
a_i\leq b_i\ \mbox{for $i=1$, \ldots, $m$.}
\]
And the element $[a_1, \ldots, a_m]$ of $\Gamma(X)$
is mapped to $\det(X_{ia_j})\in G(X)$.
It is easily verified that $\Gamma(X)$ is a 
distributive lattice.

Now we fix an $n$-dimensional vector space $V$ and a
complete flag
$V=V_n\supsetneq V_{n-1}\supsetneq\cdots\supsetneq V_1\supsetneq V_0=0$
of subspaces.
For integers $a_1$, \ldots, $a_m$ with
$1\leq a_1<\cdots< a_m\leq n$,
the Schubert subvariety $\Omega(a_1, \ldots, a_m)$
of $G_{m,n}$ is defined by
\[
\Omega(a_1, \ldots, a_m)\define
\{W\in G_{m,n}\mid
\dim(W\cap V_{a_i})\geq i\ \text{for $i=1$, \ldots, $m$}\}.
\]
If we put $b_i=n+1-a_{m+1-i}$ for $i=1$, \ldots, $m$,
$\gamma=[b_1, \ldots, b_m]$ and
$\Gamma(X;\gamma)=\{\delta\in\Gamma(X)\mid
\delta\geq\gamma\}$,
then the homogeneous coordinate ring of the Schubert variety
(Schubert cycle for short) is
\[
G(X;\gamma)\define
G(X)/(\Gamma(X)\setminus\Gamma(X;\gamma))G(X).
\]
This ring is a graded ASL over $B$ generated by $\Gamma(X;\gamma)$
(\cite{dep2}, \cite{bv}).

Set
\[
U_\gamma\define
\left(
\begin{array}{cccccccccccc}
0&\cdots&0&U_{1b_1}&\cdots&U_{1b_2-1}&U_{1b_2}&\cdots&\cdots&\cdots&\cdots&U_{1n}\\
0&\cdots&0&0&\cdots&0&U_{2b_2}&\cdots&\cdots&\cdots&\cdots&U_{2n}\\
\cdots&\cdots&\cdots&\cdots&\cdots&\cdots&\cdots&\cdots&\cdots&\cdots&\cdots\\
0&\cdots&0&0&\cdots&0&0&\cdots&0&U_{mb_m}&\cdots&U_{mn}
\end{array}
\right),
\]
where $U_{ij}$ are independent indeterminates.
Then by \cite[(7.16)]{bv}, $G(X;\gamma)\simeq G(U_\gamma)$.
So a Schubert cycle is a graded ASL generated by a
distributive lattice and an integral domain.

Now we recall the following fact.
\begin{fact}[%
%\cite{sta1}, 
\cite{sta2}, \cite{hib}]\label{fact sh}
Let $A$ be a graded ASL domain over a field generated
by a distributive lattice $D$.
Assume that
$\deg(\alpha)+\deg(\beta)=\deg(\alpha\meet\beta)+\deg(\alpha\join\beta)$
for all $\alpha$, $\beta\in D$.
Then $A$ is \gor{} if and only if
the set of all the join-irreducible elements of $D$ is pure.
\end{fact}
Recall that an element $x$ in a lattice $L$ is
said to be join-irreducible if there is unique $y$ in $L$
such that $y\covered x$,
where $y\covered x$ means that $y<x$ and there is no
element $z\in L$ such that $y<z<x$.

So in this note, we shall study the set $P$ of all the
join-irreducible elements of $\Gamma(X;\gamma)$
and give a simple proof of the criterion of the \gor{}
property of Schubert cycles.

\section{Join-irreducible elements of $\Gamma(X;\gamma)$}

As in the previous section, let $m$ and $n$ be integers
such that $1\leq m\leq n$ and
$\gamma=[b_1, \ldots, b_m]\in\Gamma(X)$.
And let $P$ be the set of all the join-irreducible elements
of $\Gamma(X;\gamma)$.

Assume that $[c_1, \ldots, c_m]$, $[d_1, \ldots, d_m]\in\Gamma(X)$.
Then it is easy to verify that
\[
[c_1,\ldots,c_m]\covered[d_1,\ldots,d_m]
\Longleftrightarrow
\exists  i;d_i=c_i+1,d_j=c_j(j\neq i).
\]
So $[c_1, \ldots, c_m]$ is a join-irreducible element of
$\Gamma(X;\gamma)$ if and only if there is a unique $i$ such that
\begin{equation}\label{join irred}
c_i>b_i,\quad c_i>c_{i-1}+1,
\end{equation}
where we set $c_0\define 0$.
For a join-irreducible element 
$[c_1, \ldots, c_m]$ of $\Gamma(X;\gamma)$, we take $i$ satisfying
(\ref{join irred}),
and set
\[
p\define n-c_i-(m-i),\quad q\define i-1.
\]
That is, $p=\#\{t\in\{1,2,\ldots,n\}\mid
t>c_i$ and $t\not\in\{c_1,c_2,\ldots,c_m\}\}$
and
$q=\#\{t\in\{1,2,\ldots,n\}\mid
t<c_i$ and $t\in\{c_1,c_2,\ldots,c_m\}\}$.
We denote the map which send $[c_1, \ldots, c_m]$
to $(p,q)$ by $\varphi$.

It is easy to construct a join-irreducible element $[c_1, \ldots, c_m]$
such that $\varphi([c_1, \ldots, c_m])=(p,q)$, if $(p,q)$
is in the image of $\varphi$.
And it also easy to verify that if
$(p,q)$ is in the image of $\varphi$
and $0\leq p'\leq p$, $0\leq q'\leq q$,
then $(p',q')$ is also in the image of $\varphi$.
Moreover,
if
$[c_1, \ldots, c_m]$ and $[d_1, \ldots, d_m]$ 
are join-irreducible elements of $\Gamma(X;\gamma)$,  and
$\varphi([c_1, \ldots, c_m])=(p,q)$, 
$\varphi([d_1, \ldots, d_m])=(p',q')$,
then
\[
[c_1, \ldots ,c_m]\leq[d_1, \ldots, d_m]\Longleftrightarrow
	p\geq p',q\geq q'.
\]
In particular, $P$ is isomorphic to a finite filter
(i.e. a subset $F$ of a poset such that 
$x\in F$ and $y>x$ implies $y\in F$)
of $\NNN\times \NNN$ whose order is defined by
\[
(p,q)\leq (p',q')\stackrel{\rm def}{\Longleftrightarrow}
p\geq p'\text{ and }q\geq q'.
\]
And if $\varphi([c_1, \ldots, c_m])=(p,q)$, then the
coheight of $[c_1, \ldots, c_m]$ is $p+q$.

Therefore, if we set
\[
\{i\mid b_i+1<b_{i+1},\ b_i<n\}=\{l_1, \ldots, l_u\},\quad
	l_1<\cdots<l_u,
\]
then the minimal elements of $P$ are
$[b_1, \ldots, b_{l_1-1},b_{l_1}+1,b_{l_1+1}, \ldots b_m]$, \ldots,
$[b_1, \ldots, b_{l_u-1}, b_{l_u}+1,b_{l_u+1},\ldots,b_m]$
and their coheights are
$n-m-b_{l_1}+2l_1-2$, \ldots,~$n-m-b_{l_u}+2l_u-2$
respectively.

\begin{example}\rm
If
$n=14$, $m=7$, $[b_1,\ldots, b_m]=[2,4,5,9,10,12,13]$, then
$u=4$, $l_1=1$, $l_2=3$, $l_3=5$ and $l_3=7$ and
the minimal elements of $P$ are
\[
\vcenter{\openup1\jot
	\halign{\hfil$\displaystyle#{}$&$\displaystyle{}#$\hfil\cr
\gamma_1=&[3,4,5,9,10,12,13]\cr
\gamma_2=&[2,4,6,9,10,12,13]\cr
\gamma_3=&[2,4,5,9,11,12,13]\cr
\gamma_4=&[2,4,5,9,10,12,14].\cr
}}
\]
And the Hasse diagram of $P$ is the following.
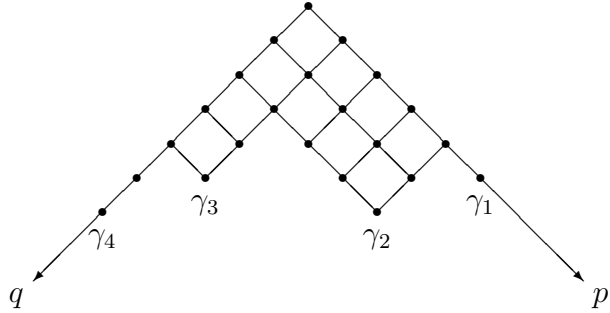
\begin{center}\unitlength\textwidth
\divide\unitlength by 300\relax
\begin{picture}(180,90)(-90,-90)
\put(0,0){\circle*{2}}
\put(-10,-10){\circle*{2}}
\put(-20,-20){\circle*{2}}
\put(-30,-30){\circle*{2}}
\put(-40,-40){\circle*{2}}
\put(-50,-50){\circle*{2}}
\put(-60,-60){\circle*{2}}
%\put(-70,-70){\circle*{2}}
\put(10,-10){\circle*{2}}
\put(0,-20){\circle*{2}}
\put(-10,-30){\circle*{2}}
\put(-20,-40){\circle*{2}}
\put(-30,-50){\circle*{2}}
\put(20,-20){\circle*{2}}
\put(10,-30){\circle*{2}}
\put(0,-40){\circle*{2}}
\put(30,-30){\circle*{2}}
\put(20,-40){\circle*{2}}
\put(10,-50){\circle*{2}}
\put(40,-40){\circle*{2}}
\put(30,-50){\circle*{2}}
\put(20,-60){\circle*{2}}
\put(50,-50){\circle*{2}}

\put(0,0){\vector(-1,-1){80}}
\put(0,0){\vector(1,-1){80}}

\put(10,-10){\line(-1,-1){40}}
\put(20,-20){\line(-1,-1){20}}
\put(30,-30){\line(-1,-1){20}}
\put(40,-40){\line(-1,-1){20}}

\put(-10,-10){\line(1,-1){40}}
\put(-20,-20){\line(1,-1){40}}
\put(-30,-30){\line(1,-1){10}}
\put(-40,-40){\line(1,-1){10}}

\put(-60,-65){\makebox(0,0)[t]{$\gamma_4$}}
\put(-30,-55){\makebox(0,0)[t]{$\gamma_3$}}
\put(20,-65){\makebox(0,0)[t]{$\gamma_2$}}
\put(50,-55){\makebox(0,0)[t]{$\gamma_1$}}

\put(-85,-85){\makebox(0,0){$q$}}
\put(85,-85){\makebox(0,0){$p$}}
\end{picture}
\end{center}
\end{example}

Note that this poset is the poset 
obtained by Stanley \cite[5.b]{sta2}.

Using Fact \ref{fact sh}, we obtain a simple proof of 
the following
criterion of \gor{} property of 
Schubert cycles.

\begin{thm}
In the notation above,
$G(X;\gamma)$ is \gor{} if and only if
$b_{l_i}-2l_i$ is constant for $i=1$, \ldots, $u$.
\end{thm}

\end{document}